\documentclass[11pt,oneside]{article}
\usepackage{amssymb,amsmath,url}
\usepackage[english]{babel}
\usepackage{pstricks,pst-node,textcomp}
\renewcommand{\le}{\leqslant}
\renewcommand{\ge}{\geqslant}
\renewcommand{\epsilon}{\varepsilon}
\renewcommand{\phi}{\varphi}

\newcommand{\RR}{\mathbb{R}}
\newcommand{\ZZ}{\mathbb{Z}}
\newcommand{\QQ}{\mathbb{Q}}

\newcommand{\NN}{\mathbb{N}}

\newcommand{\PP}{\mathbb{P}}

\newtheorem{lemma}{Lemma}
\newtheorem{theorem}{Theorem}
\newtheorem{proposition}{Proposition}
\newtheorem{definition}{Definition}

\newenvironment{proof}{\textsc{Proof. }}{\ \newline\hspace*{\fill}$\boxtimes$}

\newcommand{\bigk}{\mathop{\mathbf{K}}}

\begin{document}

\title{On effective irrationality exponents of cubic irrationals}

\author{
 Dzmitry Badziahin
}
\maketitle

\begin{abstract}
We provide an upper bound on the efficient irrationality exponents
of cubic algebraics $x$ with the minimal polynomial $x^3 - tx^2 -
a$. In particular, we show that it becomes non-trivial, i.e. better
than the classical bound of Liouville in the case $|t| > 19.71
a^{4/3}$. Moreover, under the condition $|t| > 86.58 a^{4/3}$,  we
provide an explicit lower bound on the expression $||qx||$ for all
large $q\in\ZZ$. These results are based on the recently discovered
continued fractions of cubic irrationals~\cite{badziahin_2022} and
improve the currently best-known bounds of Wakabayashi.
\end{abstract}

{\footnotesize{{\em Keywords}: cubic irrationals, continued
fractions, irrationality exponent, effective irrationality exponent

Math Subject Classification 2010: 11J68, 11J70, 11J82}}

\section{Introduction}

The irrationality exponent $\lambda(x)$ of an irrational real number
$x$ is defined as the supremum of real numbers $\lambda$ such that the
inequality
\begin{equation}\label{irr_exp}
\left|x - \frac{p}{q}\right| < q^{-\lambda}
\end{equation}
has infinitely many rational solutions $p/q$. It follows from the
classical Dirichlet theorem that for all $x\in \RR\setminus \QQ$,
$\lambda(x)\ge 2$. On the other hand, Khintchine's theorem implies
that almost all $x\in\RR$ with respect to the Lebesgue measure
satisfy $\lambda(x) = 2$. In the first half of the XX century, there
was a big interest in estimating the irrationality exponent of real
algebraic numbers. It culminated in 1955 with the work of
Roth~\cite{roth_1955}, who established the best possible result,
i.e. that for any algebraic $x\in\RR\setminus \QQ$, $\lambda(x) =
2$. Unfortunately, that result is ineffective, i.e. for $\lambda>2$
it does not allow us to find all rational $p/q$ that
satisfy~\eqref{irr_exp}. Therefore, for example, it can not be used
to solve the Thue equations
$$
F(p,q) = c
$$
in integer $p, q$, where $F\in\ZZ[x,y]$ is a homogeneous polynomial
of degree $d\ge 3$ and $c$ is some integer parameter. Since then,
many mathematicians were working on effective results regarding the
irratoinality exponents of algebraic numbers.

Given $x\in\RR\setminus \QQ$, by the effective irrationality
exponent of $x$ we define a positive real number $\lambda_{eff}(x)$
such that for all $\lambda>\lambda_{eff}(x)$ there exists an
effectively computable $Q>0$ such that all rational solutions of the
inequality~\eqref{irr_exp} in reduced form satisfy $q\le Q$.

All known upper bounds on $\lambda_{eff} (x)$ for algebraic real $x$
are much weaker than in Roth's theorem. First of all, the classical
theorem of Liouville states that $\lambda_{eff}(x) \le d$ where $d$
is the degree of $x$. Therefore any non-trivial bound on
$\lambda_{eff}(x)$ should be strictly smaller than $d$. One of the
notable improvements of Liouville's bound is based on Feldman's
refinement of the theory of linear forms in
logarithms~\cite{feldman_1968}. Its advantage is that it gives
$\lambda_{eff}(x)<d$ for all algebraic numbers. However, the
difference between $\lambda_{eff}(x)$ and $d$ is usually extremely
small. For state-of-the-art results regarding this approach, we
refer to the book of Bugeaud~\cite{bugeaud_2018}. For other notable
achievements on this problem, we refer to~\cite{bennett_1997,
bpv_1996} and the references therein.

In this paper, we focus on the case of cubic irrationals.
Multiplying by some integer number and shifting by another rational
number, we can always guarantee that the minimal polynomial of the
resulting cubic $x$ is $x^3 + px + q$ for some $p,q\in\ZZ$. Notice
also that such a transformation does not change the (effective)
irrationality exponent of $x$. The first general result about
$\lambda_{eff}(x)$ for these specific values $x$ was achieved by
Bombieri, van~der~Poorten and Vaaler~\cite{bpv_1996} in 1996. They
showed that under the conditions $|p|>e^{1000}$ and $|p|\ge q^2$,
one has
$$
\lambda_{eff}(x)\le \frac{2\log(|p|^3)|}{\log (|p|^3/q^2)} +
\frac{14}{(\log (|p|^3/q^2))^{1/3}}.
$$
Later, Wakabayashi~\cite{wakabayashi_2002} improved that bound and
showed that $\lambda_{eff}(x)\le \lambda_w(p,q)$. It becomes
non-trivial (i.e. smaller than 3) under the condition
\begin{equation}\label{wak_ineq}
|p|>2^{2/3}3^4|q|^{8/3} \left(1+\frac{1}{390|q|^3}\right)^{2/3}
\end{equation}
and for large $p$ and $q$ it asymptotically behaves like
$\lambda_w(p,q)\sim 2 + (4\log|q| + 2\log 108)/(3\log |p|)$.

In this paper we investigate what estimates on $\lambda_{eff}(x)$
can be achieved with help of the convergents of the recently
discovered continued fractions~\cite{badziahin_2022} of cubic
irrationals with the minimal polynomial $x^3 - tx^2 - a\in\ZZ[x]$.
It was shown there that, as soon as $|t|^3
> 12 a>0$, the real root of this equation with the largest absolute
value admits the continued fraction {\footnotesize
$$
x = \bigk\left[ \begin{array}{l}  \\ \!\!t \end{array} \overline{\begin{array}{llll}3(12k+1)(3k+1)\alpha&3(12k+5)(3k+2)\alpha&3(12k+7)(6k+5)\alpha& 3(12k-1)(6k+1)\alpha\\
(2i+1)t^2&(2i+1)t&2(2i+1)t^2&(2i+1)t\!\!
\end{array}\hspace{-1.5ex}}\;\; \right].
$$} Here $i$ is the index of the corresponding partial quotient and
$k = \big[\frac{i}{4}\big]$. Notice that the change of variables $y
= \frac{q}{x}$ transforms cubic $x$ from~\cite{bpv_1996,
wakabayashi_2002} to numbers in this paper with $t = -p$ and $a =
-q^2$.

The resulting upper bounds on $\lambda_{eff}(x)$ (see
Theorems~\ref{th1} and~\ref{th2}) depend on prime factorisations of
$a$ and $t$ but in any case they are better than those
in~\cite{wakabayashi_2002}. Theorem~\ref{th1} states that under the
condition $|t|^3>12|a|$ the largest real root of the above cubic
equation satisfies
$$
||qx|| \ge \tau(t,a) q^{-\lambda(t,a)} (\log
(8|t_1|q))^{\lambda(t,a)-1/2},\quad \mbox{for all } q>Q_0(t,a),
$$
where all values for $\tau(t,a), \lambda(t,a)$ and $q_0(t,a)$ are
explicitly provided. This inequality always gives a non-trivial
upper bound on $\lambda_{eff}(x)$ under the condition $|t|> 86.58
a^{4/3}$. It translates to the condition $|p| > 86.58 |q|^{8/3}$
which is better than~\eqref{wak_ineq}, where even without the term
in brackets we have $|p|>128.57 |q|^{8/3}$. Moreover, the parameter
$a$ does not have to be of the form $-q^2$, therefore many cubics
from this paper do not plainly transfer to those in~\cite{bpv_1996,
wakabayashi_2002}. On top of that, Theorem~\ref{th2} provides an
even better upper bound on $\lambda_{eff}(x)$ but does not give an
explicit lower bound on $||qx||$ as above. That bound always becomes
non-trivial as soon as $|t|>19.71 a^{4/3}$. These bounds between $t$
and $a$ are obtained in Section~\ref{sec3}.

\section{Preliminaries and main results}

Consider the real root of the equation $x^3 - tx^2 - a=0$ that has
the largest absolute value among other roots of the same equation.
Notice that $a\neq 0$ because otherwise $x$ is not cubic. Also, by
replacing $x$ with $-x$ if needed, we can guarantee that $a>0$.

Next, by the standard CF transformations (see \cite[Lemma
1]{badziahin_2022}), we can cancel some common divisors of $t$ and
$3a$ from the continued fraction of $x$. Let $g_1:= \gcd(t^2, 3a)$
and $g_2:=\gcd(t, 3a/g_1)$. For convenience, we denote $t^2 =
g_1t_2$, $t = g_2t_1$ and $3a = g_1g_2a^*$. We divide the following
partial quotients of $x$ by $g_1$: $\beta_1,a_1,\beta_2$; $\beta_3,
a_3,\beta_4$; \ldots, $\beta_{2k-1},a_{2k-1},\beta_{2k}$,\ldots.
After that we divide the following partial quotients by $g_2$:
$\beta_0, a_0, \beta_1$; $\beta_2,a_2,\beta_3$; \ldots; $\beta_{2k},
a_{2k}, \beta_{2k+1}$,\ldots The resulting continued fraction has
the same limit $x$ as the initial one. To make $\beta_0$ integer, we
consider the number $x/g_2$ instead of $x$. Its continued fraction
is then {\footnotesize
\begin{equation}\label{no4}
\bigk\left[ \begin{array}{l}  \\ \!\!t_1 \end{array} \begin{array}{llll}(12k+1)(3k+1)a^*&(12k+5)(3k+2)a^*&(12k+7)(6k+5)a^*& (12k-1)(6k+1)a^*\\
(2i+1)t_2&(2i+1)t_1&2(2i+1)t_2&(2i+1)t_1\!\!
\end{array}\hspace{-1.5ex}\;\; \cdots\right]
\end{equation}}

We define the following notions:
\begin{equation}\label{def_c6}
c_6 = c_6(t_1,t_2,a^*):= \left\{ \begin{array}{ll}
\displaystyle\frac{\sqrt{64t_1t_2 + 270a^*}}{2^{7/4}ec_1}&\mbox{if }\; t>0\\[3ex]
\displaystyle\frac{\sqrt{64|t_1t_2| - 54a^*}}{2^{7/4}ec_1}&\mbox{if
}\; t<0,
\end{array}\right.
\end{equation}

\begin{equation}\label{def_c7}
c_7 = c_7(t_1,t_2,a^*):= \left\{ \begin{array}{ll}
\displaystyle\frac{2^{7/4}ec_1(16t_1t_2+9a^*)}{9a^*\sqrt{64t_1t_2 + 270a^*}}&\mbox{if }\; t>0\\[3ex]
\displaystyle\frac{2^{23/4}ec_1(|t_1t_2|-3a^*)}{9a^*\sqrt{64|t_1t_2|
- 54a^*}}&\mbox{if }\; t<0,
\end{array}\right.
\end{equation}
where $c_1 = c_1(a^*)$ is defined in~\eqref{def_c1}. Next,
\begin{equation}\label{def_tau}
\tau = \tau(t_1,t_2,a^*):= \frac{g_2\cdot(\log
c_7)^{1/2}}{8c_6^8(8|t_1|)^{\frac{\log c_6}{\log c_7}}};\quad q_0 =
q_0(t_1,t_2,a^*) := \frac{c_7^4}{8|t_1|}.
\end{equation}

The main result of this paper is

\begin{theorem}\label{th1}
For all integer $q\ge q_0$ one has
$$
||qx|| > \tau q^{-\lambda} \log (8|t_1| q)^{-\lambda-1/2}
$$
where $\lambda = \frac{\log c_6}{\log c_7}$. In particular,
$\lambda_{eff}(x) \le \lambda$.
\end{theorem}

As shown in Section~\ref{sec4}, the constant $c_1$ can be replaced
by a bigger constant $c_2$ defined in~\eqref{def_c2}. But in that
case, writing such an explicit inequality as in Theorem~\ref{th1} is
much harder (but theoretically possible). We do not provide it in
the exact form here but only state the following result. Let $c_6^*$
and $c_7^*$ be defined in the same way as $c_6$ and $c_7$ but with
the constant $c_2$ instead of $c_1$.

\begin{theorem}\label{th2}
The effective irrationality exponent of $x$ satisfies
$$
\lambda_{eff} (x)\le \lambda^* := \frac{\log c_6^*}{\log c_7^*}.
$$
\end{theorem}

\section{Analysis of the results}\label{sec3}

Theorem~\ref{th1} provides a nontrivial lower bound for $||qx||$ as
soon as $\frac{\log c_6}{\log c_7}$ is strictly less that 2 or in
other words, $c_6<c_7^2$. In view of~\eqref{def_c6}
and~\eqref{def_c7}, for $t>0$ this is equivalent to
$$
\frac{\sqrt{64t_1t_2+270a^*}}{2^{7/4} ec_1} <
\frac{2^{7/2}e^2c_1^2(16t_1t_2+9a^*)^2}{81a^{*2}(64t_1t_2+270a^*)}
$$
$$
\Longleftrightarrow\qquad a^{*2}(64t_1t_2 + 270a^*)^{3/2} <
\frac{2^{21/4}e^3c_1^3}{81}(16t_1t_2+9a^*)^2.
$$
Define the parameter $u$ such that $16t_1t_2 = ua^{*4}-9a^*$. Also
for convenience define $\tau:= \frac{2^{21/4}e^3c_1^3}{81}$. Then
the last inequality rewrites
$$
a^{*2}(4ua^{*4} + 234 a^*)^{3/2} < \tau u^2a^{*8} \quad
\Longleftrightarrow\quad 4ua^{*3} + 234 < \tau^{2/3}u^{4/3} a^{*3}.
$$
Notice that for $u\ge \left(\frac{4}{\tau^{2/3}} +
\frac{234}{64a^{*3}}\tau^{4/3}\right)^3$ one has
$$
\tau^{2/3}u^{4/3}a^{*3}\ge \left(4 +
\frac{234}{64a^{*3}}\tau^2\right)\cdot ua^{*3} \ge 4ua^{*3} +
\frac{234\tau^2}{64a^{*3}}\cdot \frac{4^3}{\tau^2}a^{*3} = 4ua^{*3}
+ 234
$$
and the condition $c_6 < c_7^2$ is satisfied.

Recall that $t_1t_2 = t^3/(g_1g_2)$ and $a^* = 3a/(g_1g_2)$.
Therefore, $16t_1t_2 > ua^{*4} - 9a^*$ is equivalent to $16t^3 >
\frac{81}{(g_1g_2)^3}ua^4-27a$.

From the definition~\eqref{def_c1} we see that $c_1$ and hence
$\tau$ depends on the prime factorisation of $a^*$. If $3\nmid a^*$
then we always get $c_1> 0.0924$ which in turn implies $\tau >
0.00744$ and
$$
\left(\frac{4}{\tau^{2/3}} +
\frac{234}{64a^{*3}}\tau^{4/3}\right)^3\le 104.97^3.
$$
On the other hand, we have $g_1g_2\ge 3$.

Finally, we get that for $t>0$ the non-trivial bound on
$\lambda_{eff}(x)$ is always achieved if $t >
\frac{104.97}{3}\cdot\left(\frac{81}{16}\right)^{1/3}a^{4/3}\approx
60.08 a^{4/3}$, however for many pairs $a$ and $t$ it is satisfied
under essentially weaker conditions.

In the case $3\mid a^*$ we get $c_1> 0.13329$, thus $\tau > 0.0223$
and
$$
\left(\frac{4}{\tau^{2/3}} +
\frac{234}{64a^{*3}}\tau^{4/3}\right)^3\le 50.42^3.
$$
Then the non-trivial bound is achieved if $t>
50.42\cdot\left(\frac{81}{16}\right)^{1/3}a^{4/3}\approx 86.57
a^{4/3}$.

The case $t<0$ is dealt analogously. The condition $c_6<c_7^2$ is
equivalent to
$$
\frac{\sqrt{64|t_1t_2|-54a^*}}{2^{7/4}ec_1} <
\frac{2^{7/2}e^2c_1^2(16|t_1t_2|-48a^*)^2}{81a^{*2}(64|t_1t_2|-54a^*)}
$$
$$
\Longleftrightarrow\qquad a^{*2}(64|t_1t_2| - 54a^*)^{3/2} <
\tau(16|t_1t_2|-48a^*)^2.
$$
Define $u$ such that $16|t_1t_2| = ua^{*4} + 48a^*$. Then the last
inequality rewrites
$$
4ua^{*3} + 138 < \tau^{2/3}u^{4/3}a^{*3}.
$$
One can check that the last inequality is satisfied for $u\ge
\left(\frac{4}{\tau^{2/3}} + \frac{138
\tau^{4/3}}{64a^{*3}}\right)^3$. As in the case of positive $t$, the
right hand side is always smaller than $104.93^3$ in the case
$3\nmid a^*$ and is smaller than $50.42^3$ in the case $3\mid a^*$.
Therefore for $3\nmid a^*$ the condition $c_6<c_7^2$ is always
satisfied if $|t|^3
> \frac{104.93^3\cdot 81}{3^3\cdot 16} a^4 + 9a$ which follows from $|t|
>
60.06 a^{4/3}$. For $3\mid a^*$, similar computations give us $|t|^3
> \frac{50.42^3\cdot 81}{16} a^4 + 9a$ which follows from $|t|>
86.58a^{4/3}$.

One can repeat the same analysis as above for Theorem~\ref{th2}. In
that case, the constant $c_1$ in the computations should be replaced
by $c_2$. We have that if $3\nmid a^*$ it always satisfies $c_2\ge
0.1939$, which in turn implies $\tau\ge 0.0688$ and then
$$
\left(\frac{4}{\tau^{2/3}} +
\frac{234}{64a^{*3}}\tau^{4/3}\right)^3\le 23.93^3.
$$
If $3\mid a^*$, we have $c_2\ge 0.2797$, $\tau \ge 0.206$ and
$$
\left(\frac{4}{\tau^{2/3}} +
\frac{234}{64a^{*3}}\tau^{4/3}\right)^3\le 11.47^3.
$$

Finally, for the case $3\nmid a^*$, the condition $c_6<c_7^2$ is
always satisfied if $|t|^3
> \frac{23.93^3\cdot 81}{3^3\cdot 16}a^4 + 9a$ which follows from
$|t|>13.72 a^{4/3}$. For the case $3\mid a^*$ similar calculations
give $|t|> 19.71 a^{4/3}$.

\section{Nice, convenient and perfect continued fractions}

\begin{definition} Let $x$ be a continued fraction given by
$$
x = \bigk \left[
\begin{array}{ccc}
\beta_0&\beta_1&\beta_2\\
a_0&a_1&a_2
\end{array} \cdots
\right] ;\quad \beta_i, a_i\in\ZZ,\; \forall i\in\ZZ_{\ge 0}.
$$
For given positive integers $k,r$ with $-1\le r\le k$ we define
$$
\gamma_{k,r}:= \prod_{i=k-r\atop 2\mid (i-k+r)}^{k+r} \beta_i, \quad
\gamma_{k,-1}:=1.
$$
We say that $x$ is \textbf{$d$-nice} at index $k$ where $d,k\in\NN$
if $d\mid a_k$ and for all positive integer $r\le k$ one has
$a_{k-r}\beta_{k+r}\gamma_{k,r-2}\equiv -a_{k+r}\gamma_{k,r-1}$ (mod
$d$). We call $x$ \textbf{$(d,r)$-perfect} at index $k$, where $0\le
r\le k$ if it is $d$-nice at index $k$ and
$\beta_{k-r}\equiv\beta_{k+r+1}\equiv 0$ (mod $d$).
\end{definition}

We say that $x$ is \textbf{eventually $d$-nice} at index $k$ from
position $k_0$ if the same conditions as above are satisfied for all
$0\le r\le k-k_0$. In this paper the value $k_0$ will often be a
fixed absolute constant. If there is no confusion about its value we
will omit it in the text.

Let $p_n/q_n$ be the $n$'th convergent of $x$. Define the following
matrices
$$
S_n:= \left(\begin{array}{cc} p_{n}&q_{n}\\
p_{n-1}&q_{n-1}
\end{array}\right);\qquad C_n:=\left(\begin{array}{cc} a_{n}&\beta_{n}\\
1&0
\end{array}\right).
$$
From the theory of continued fractions we know that $S_n =
C_nC_{n-1}\cdots C_0$. To make this product shorter, we use the
usual notation but in the descending order: $ S_n = \prod_{i=n}^0
C_i$.

\begin{lemma}\label{lem1}
Let the continued fraction $x$ be eventually $d$-nice at index $k$
from the position $k_0$. Then for all $0\le r\le k-k_0$ one has
$$
\prod_{i=k+r}^{k-r} C_i \equiv \left(\begin{array}{cc}
0& \gamma_{k,r}\\
\gamma_{k,r-1}&0
\end{array}\right)\;\;(\mathrm{mod}\; d).
$$
Moreover, if $x$ is $(d,r)$-perfect at index $k$ then
$\prod_{i=k+r+1}^{k-r} C_i  \equiv \mathbf{0}$ (mod $d$).
\end{lemma}

\proof

We prove by  induction on $r$. For $r=0$ the statement is
straightforward. Suppose that the statement is true for $r$ and
verify it for $r+1$.
$$
\begin{array}{rl}
\displaystyle\prod_{i=k+r+1}^{k-r-1} C_i &\equiv
\left(\begin{array}{cc}
a_{k+r+1}&\beta_{k+r+1}\\
1&0
\end{array}\right)
\left(\begin{array}{cc}
0& \gamma_{k,r}\\
\gamma_{k,r-1}&0
\end{array}\right)  \left(\begin{array}{cc}
a_{k-r-1}&\beta_{k-r-1}\\
1&0
\end{array}\right)\\
&\equiv \left(\begin{array}{cc}
a_{k-r-1}\beta_{k+r+1}\gamma_{k,r-1} + a_{k+r+1}\gamma_{k,r}& \beta_{k-r-1}\beta_{k+r+1}\gamma_{k,r-1}\\
\gamma_{k,r}&0
\end{array}\right)\;\;(\mathrm{mod}\; d).
\end{array}
$$
By the conditions of $d$-nice CF at index $k$, the last matrix is congruent
to $\left(\begin{array}{cc}
0& \gamma_{k,r+1}\\
\gamma_{k,r}&0
\end{array}\right)$.

If $x$ is $(d,r)$-perfect at index $k$ then $\gamma_{k,r}\equiv 0$
(mod $d$) and we have
$$
\prod_{i=k+r+1}^{k-r} C_i \equiv \left(\begin{array}{cc}
a_{k+r+1}&0\\
1&0
\end{array}\right)
\left(\begin{array}{cc}
0& \gamma_{k,r}\\
\gamma_{k,r-1}&0
\end{array}\right)\equiv \mathbf{0}\;\;(\mathrm{mod}\; d).
$$
\endproof

Another two properties of $d$-nice continued fractions that easily follow
from the definition are
\begin{itemize}
\item Let $d_1,d_2$ be two coprime positive integer numbers. If a continued
    fraction is eventually $d_1$-nice and eventually $d_2$-nice at the same
    index $k$ for the same position $k_0$ then it is eventually
    $d_1d_2$-nice at index $k$.
\item If a continued fraction is eventually $d$-nice at index $k$ then it is also
    eventually $e$-nice at the same index from the same position for all positive
    integer divisors $e$ of $d$.
\end{itemize}

\begin{definition}
We say that the continued fraction $x$ is (eventually)
$d$-convenient at index $k$ if there exists a sequence $(c_r)_{0\le
r\le \lfloor k/2\rfloor}$ of residues modulo $m$ such that for all
positive integers $r\le k$ (resp. $r\le k-k_0$) one has
\begin{itemize}
\item $\beta_{k+r+1}\equiv
c_{\lfloor \frac{r}{2}\rfloor}\beta_{k-r}$ (mod $d$);
\item if $r$ is odd then $a_{k+r}\equiv -c_{\lfloor \frac{r}{2}\rfloor} a_{k-r}$ (mod
$d$);
\item if $r$ is even then $a_{k+r}\equiv -a_{k-r}$ (mod $d$).
\end{itemize}
\end{definition}

\begin{lemma}\label{lem3}
Let $d>2$. Then any eventually $d$-convenient continued fraction at
index $k$ is also eventually $d$-nice. For $d=2$, any $d$-convenient
continued fraction at index $k$ such that $a_k\equiv 0$ (mod $d$) is
also $d$-nice.
\end{lemma}
\proof First of all, for $d>2$ and $r=0$ the condition
$a_{k+r}\equiv -a_{k-r}$ (mod $d$) implies that $a_k\equiv 0$ (mod
$d$), which is the first condition of $d$-nice CF.

Secondly, one can check that the first condition of $d$-convenient
CF implies that for odd $r$, $\gamma_{k,r} \equiv c_{\lfloor
\frac{r}{2}\rfloor}\gamma_{k,r-1}\beta_{k-r}\equiv
\gamma_{k,r-1}\beta_{k+r+1}$ (mod $d$). Then we get
$$
a_{k-r-1}\beta_{k+r+1}\gamma_{k,r-1}\equiv
a_{k-r-1}\gamma_{k,r}\equiv -a_{k+r+1}\gamma_{k,r}\;(\mathrm{mod}\;
d)
$$
and the second condition of $d$-nice CF is verified.

Thirdly, for even $r$ we get $c_{r/2}\gamma_{k,r}\equiv
c_{r/2}\gamma_{k,r-1}\beta_{k-r}\equiv \gamma_{k,r-1}\beta_{k+r+1}$
(mod $d$) and therefore
$$
a_{k-r-1}\beta_{k+r+1}\gamma_{k,r-1}\equiv
c_{r/2}a_{k-r-1}\gamma_{k,r}\equiv
-a_{k+r+1}\gamma_{k,r}\;(\mathrm{mod}\; d).
$$
Again, the second condition of $d$-nice CF is satisfied.
\endproof

\section{Divisibility of entries of $S_n$}\label{sec4}

\begin{lemma}\label{lem2}
Let $k\in \NN$, $k\ge 2$ and $d$ be any integer divisor of \ $2k+1$. The
continued fraction~\eqref{no4} is eventually $d$-convenient at index $k$ from
the position 2. Additionally, the same statement is true for $k\equiv 3$ (mod
4) and $d=2$.
\end{lemma}

In the further discussion we will always deal with eventually
$d$-convenient or $d$-nice continued fractions from the position 2.
Therefore, to make the description shorter, we will omit the words
`eventually' and `from the position 2' and call the continued
fraction~\eqref{no4} $d$-convenient or $d$-nice.

\proof We will check the conditions of $d$-convenient continued fraction
separately for each of the cases, depending on $k$ modulo 4.

{\bf Case $k=4m+1$.} Then $m\equiv -\frac38$ (mod $d$) and we
use~\eqref{no4} to compute
$$
a_{k+4r}\equiv 8rt_2,\quad \beta_{k+4r}=
a^*(12m+1+12r)(3m+1+3r)\equiv
\frac{a^*}{16}(24r-7)(24r-1)\;(\mathrm{mod}\; d);
$$$$
a_{k+4r+1}\equiv (8r+2)t_1,\quad \beta_{k+4r+1}\equiv
\frac{a^*}{16}(24r+1)(24r+7)\;(\mathrm{mod}\; d);
$$$$
a_{k+4r+2}\equiv 2(8r+4)t_2,\quad \beta_{k+4r+2}\equiv
\frac{a^*}{8}(24r+5)(24r+11)\;(\mathrm{mod}\; d);
$$$$
a_{k+4r-1}\equiv (8r-2)t_1,\quad \beta_{k+4r-1}\equiv
\frac{a^*}{8}(24r-11)(24r-5)\;(\mathrm{mod}\; d).
$$
The conditions of $d$-convenient continued fraction at index $k$ can
now be easily checked where $c_r$ is the constant 1 sequence.

We proceed the same way in all other cases.

{\bf Case $k=4m+2$.} Then $m\equiv -\frac{5}{8}$ (mod $d$) and
$$
a_{k+4r}\equiv 8rt_1,\quad \beta_{k+4r}\equiv
\frac{a^*}{16}(24r-5)(24r+1)\;(\mathrm{mod}\; d);
$$$$
a_{k+4r+1}\equiv 2(8r+2)t_2,\quad \beta_{k+4r+1}\equiv
\frac{a^*}{8}(24r-1)(24r+5)\;(\mathrm{mod}\; d);
$$$$
a_{k+4r+2}\equiv (8r+4)t_1,\quad \beta_{k+4r+2}\equiv
\frac{a^*}{8}(24r+7)(24r+13)\;(\mathrm{mod}\; d);
$$$$
a_{k+4r-1}\equiv (8r-2)t_2,\quad \beta_{k+4r-1}\equiv
\frac{a^*}{16}(24r-13)(24r-7)\;(\mathrm{mod}\; d).
$$
One can easily check that for $s\equiv 0,1$ (mod 4), $\beta_{k+s+1}
\equiv 2\beta_{k-s}$ (mod $d$) and for $s\equiv 2,3$ (mod 4),
$\beta_{k+s+1}\equiv 2^{-1} \beta_{k-s}$ (mod $d$). Also, for
$s\equiv 1$ (mod 4), $a_{k+s} \equiv 2a_{k-s}$ (mod $d$) and for
$s\equiv 3$ (mod 4), $a_{k+s} \equiv 2^{-1}a_{k-s}$ (mod $d$).
Hence, the conditions of $d$-convenient CF are verified, where the
sequence $c_s$ is periodic with the period $2, 2^{-1}$.

{\bf Case $k=4m+3$, $d\neq 2$.} Then $m\equiv -\frac{7}{8}$ (mod
$d$) and
$$
a_{k+4r}\equiv 16rt_2,\quad \beta_{k+4r}\equiv
\frac{a^*}{8}(24r-7)(24r-1)\;(\mathrm{mod}\; d);
$$$$
a_{k+4r+1}\equiv (8r+2)t_1,\quad \beta_{k+4r+1}\equiv
\frac{a^*}{8}(24r+1)(24r+7)\;(\mathrm{mod}\; d);
$$$$
a_{k+4r+2}\equiv (8r+4)t_2,\quad \beta_{k+4r+2}\equiv
\frac{a^*}{16}(24r+5)(24r+11)\;(\mathrm{mod}\; d);
$$$$
a_{k+4r-1}\equiv (8r-2)t_1,\quad \beta_{k+4r-1}\equiv
\frac{a^*}{16}(24r-11)(24r-5)\;(\mathrm{mod}\; d).
$$
One can then check the conditions of $d$-convenient CF at index $k$
for the constant 1 sequence $c_r$.

{\bf Case $k=4m+3$, $d=2$.} in this case one can easily see that
$a_{k+4r}\equiv 0$ (mod 2), $a_{k+4r+1}\equiv a_{k+4r+3}\equiv t_1$
(mod 2), $a_{k+4r+2}\equiv t_2$ (mod 2); $\beta_{k+4r}\equiv
\beta_{k+4r+1}\equiv a^*$ (mod 2) and $\beta_{k+4r+2}\equiv
\beta_{k+4r-1}$ (mod 2). Therefore, the CF is 2-convenient at index
$k$ with the constant 1 sequence $c_r$.

{\bf Case $k=4m$.} Then $m\equiv -\frac{1}{8}$ (mod $d$) and
$$
a_{k+4r}\equiv 8rt_1,\quad \beta_{k+4r}\equiv
\frac{a^*}{8}(24r-5)(24r+1)\;(\mathrm{mod}\; d);
$$$$
a_{k+4r+1}\equiv (8r+2)t_2,\quad \beta_{k+4r+1}\equiv
\frac{a^*}{16}(24r-1)(24r+5)\;(\mathrm{mod}\; d);
$$$$
a_{k+4r+2}\equiv (8r+4)t_1,\quad \beta_{k+4r+2}\equiv
\frac{a^*}{16}(24r+7)(24r+13)\;(\mathrm{mod}\; d);
$$$$
a_{k+4r-1}\equiv 2(8r-2)t_2,\quad \beta_{k+4r-1}\equiv
\frac{a^*}{8}(24r-13)(24r-7)\;(\mathrm{mod}\; d).
$$
One can then check the conditions of $d$-convenient CF with the
periodic sequence $c_r$ with the period $2^{-1},2$. \endproof

Lemmata~\ref{lem3} and~\ref{lem2} show that $x$ is $d$-nice at each
index $k\ge 2$ for appropriately chosen $d$. As the next step, we
show that for almost every prime $p$, it is also $(p,t)$-perfect at
infinitely many carefully chosen indices $k$ and $t$. This fact will
allow us to show that all the entries of $C_kC_{k-1}\cdots C_1$ are
multiples of some big power of $p$.

First of all, let's consider the case $p>2$ and $p\mid a^*$. Let
$s\in\NN$ be such that $p^s|| a^*$. Consider $q = p^l$ for some
$1\le l\le s$. If we write $q=2m+1$ then we have $q\mid \alpha_k$
for $k=m + rq = \frac{(2r+1)q - 1}{2}$ where $r\in\ZZ_{\ge 0}$. One
can easily see that for any such value of $k$, $x$ is
$(q,0)$-perfect at index $k$. In view of Lemma~\ref{lem1}, we can
then split the product $S_k$ into $\left\lfloor
\frac{2k+q-1}{2q}\right\rfloor$ groups such that all entries of the
resulting product matrix in each group are multiples of $q$.
Finally, we combine this information for $q^l$ for all $1\le l\le s$
and derive that all entries of the product $\prod_{i=k}^2 C_i$ are
divisible by
$$
p^{\;\sum\limits_{i=1}^s \left\lfloor
\frac{2k+p^i-1}{2p^i}\right\rfloor}.
$$

Next, consider the case $p=2$ and $p\mid a^*$. We have $p\mid a_k$
for all $k\equiv 3$ (mod 4) and one can easily see that for all such
$k$, $x$ is $(p,0)$ perfect at index $k$. Then the analogous
application of Lemma~\ref{lem1} as in the previous case implies that
all entries of $\prod_{i=k}^2 C_i$ are divisible by $ 2^{\lfloor
k/4\rfloor }.$

For the case $p=2$, $p\nmid a^*$ the result is slightly weaker.
From~\eqref{no4} one can verify that $\beta_{8m+2}\equiv
\beta_{8m+5}\equiv 0$ (mod 2) for all $m\in\ZZ_{\ge 0}$ and
therefore $x$ is $(2,1)$-perfect at indices $8m+3$. Then Lemma
\ref{lem1} then implies that all entries of $\prod_{i=k}^2 C_i$ are
divisible by $2^{\lfloor (k+3)/8\rfloor }.$

Finally, in the next lemma we consider the remaining case of $p\in\NN$ that
do not divide $a^*$.

\begin{lemma}\label{lem4}
Let $p\in\NN$ be such that $\gcd(p,6)=1$. Then for all $k\in\ZZ$ all
the entries of the product of matrices $\prod_{i=k}^2 C_i$ are
divisible by
$$
p^{\left\lfloor \frac{3k + p-2}{3p}\right\rfloor}.
$$
\end{lemma}

\proof We prove by routinely considering all the cases, depending on $p$
modulo 12.

{\bf Case $p = 12m+1$.} Then with help of~\eqref{no4} one can verify that for
all $r\in\ZZ_{\ge 0}$,
$$
\begin{array}{rl}
0&\equiv \beta_{4(m+rp)+1}\equiv
\beta_{4(2m+rp)}\equiv\beta_{4(4m+rp)+1}\equiv \beta_{4(5m+rp)+2}\equiv
\beta_{4(7m+rp)+3}\equiv \beta_{4(8m+rp)+2}\\[1ex]
&\equiv \beta_{4(10m+rp)+3}\equiv \beta_{4(11m+rp)+4}\;(\mathrm{mod}\; p)
\end{array}
$$
and $0\equiv a_{6m+rp}$ (mod $p$). In view of Lemma~\ref{lem2}, we then
derive that $x$ is $(p,2m-1)$-perfect at indices $k = 6m+2rp$ for all
$r\in\ZZ_{\ge 0}$ and is $(p,2m)$-perfect at indices $k = 6m+(2r+1)p$.
Lemma~1 then implies that all the entries of the following products of
matrices are divisible by $p$:
$$
\prod_{i=8m+2rp}^{4m+2rp+1} C_i,\qquad \prod_{i=20m+2rp+2}^{16m+4rp+1} C_i.
$$
Finally, one can easily check that for $k = (n+1)p - \frac{p-1}{3}$,
the product $\prod_{i=k}^2 C_i$ contains $n+1$ blocks of the above
form. Therefore all its entries are divisible by $p^{n+1}$.

The other cases are done analogously.

{\bf Case $p=12m+5$.} One verifies that for all $r\in\ZZ_{\ge 0}$,
$$
\begin{array}{rl}
0&\equiv \beta_{4(m+rp)+2}\equiv
\beta_{4(2m+rp)+3}\equiv\beta_{4(4m+rp)+6}\equiv \beta_{4(5m+rp)+9}\equiv
\beta_{4(7m+rp)+12}\equiv \beta_{4(8m+rp)+13}\\[1ex]
&\equiv \beta_{4(10m+rp)+16}\equiv \beta_{4(11m+rp)+19}\;(\mathrm{mod}\; p)
\end{array}
$$
and $0\equiv a_{6m+rp+2}$ (mod $p$). Then Lemma~\ref{lem2} implies
that $x$ is $(p,2m)$-perfect at indices $k=6m+2rp$ for all $r\in\ZZ$
and is $(p,2m-1)$-perfect at indices $k=6m+(2r+1)p$.
Lemma~\ref{lem1} then implies that all the entries of the following
products are divisible by $p$:
$$
\prod_{i=8m+2rp+3}^{4m+2rp+2} C_i,\qquad \prod_{i=20m+2rp+9}^{16m+4rp+6} C_i.
$$
For $k\ge (n+1)p-\frac{p-2}{3}$ one can easily check that the product
$\prod_{i=k}^2 C_i$ contains $n+1$ blocks of the above form. Therefore all
its entries are divisible by $p^{n+1}$.

{\bf Case $p=12m+7$.} Then for all $r\in\ZZ_{\ge 0}$,
$$
\begin{array}{rl}
0&\equiv \beta_{4(m+rp)+3}\equiv
\beta_{4(2m+rp)+4}\equiv\beta_{4(4m+rp)+9}\equiv \beta_{4(5m+rp)+12}\equiv
\beta_{4(7m+rp)+17}\equiv \beta_{4(8m+rp)+18}\\[1ex]
&\equiv \beta_{4(10m+rp)+23}\equiv \beta_{4(11m+rp)+26}\;(\mathrm{mod}\; p)
\end{array}
$$
and $0\equiv a_{6m+rp+3}$ (mod $p$). Lemmata~\ref{lem2} and~\ref{lem1} imply
that all the entries of the following products are divisible by $p$:
$$
\prod_{i=8m+2rp+4}^{4m+2rp+3} C_i,\qquad \prod_{i=20m+2rp+12}^{16m+4rp+9} C_i.
$$
For $k\ge (n+1)p-\frac{p-1}{3}$ one can easily check that the product
$\prod_{i=k}^2 C_i$ contains $n+1$ blocks of the above form. Therefore all
its entries are divisible by $p^{n+1}$.

{\bf Case $p=12m+11$.} Then for all $r\in\ZZ_{\ge 0}$,
$$
\begin{array}{rl}
0&\equiv \beta_{4(m+rp)+4}\equiv
\beta_{4(2m+rp)+7}\equiv\beta_{4(4m+rp)+14}\equiv \beta_{4(5m+rp)+19}\equiv
\beta_{4(7m+rp)+26}\equiv \beta_{4(8m+rp)+29}\\[1ex]
&\equiv \beta_{4(10m+rp)+36}\equiv \beta_{4(11m+rp)+41}\;(\mathrm{mod}\; p)
\end{array}
$$
and $0\equiv a_{6m+rp+5}$ (mod $p$). Lemmata~\ref{lem2} and~\ref{lem1} imply
that all the entries of the following products are divisible by $p$:
$$
\prod_{i=8m+2rp+7}^{4m+2rp+4} C_i,\qquad \prod_{i=20m+2rp+19}^{16m+4rp+14} C_i.
$$
For $k\ge (n+1)p-\frac{p-2}{3}$ one can easily check that the product
$\prod_{i=k}^2 C_i$ contains $n+1$ blocks of the above form. Therefore all
its entries are divisible by $p^{n+1}$.

In all four cases we have that for $k\ge (n+1)p - \frac{p-2}{3}$ all
the entries of $\prod_{i=k}^2 C_i$ are divisible by $p^{n+1}$.
Writing it in terms of $k$ we get that this power of $p$ is
$$
\left\lfloor \frac{k+\frac{p-2}{3}}{p}\right\rfloor = \left\lfloor
\frac{3k + p - 2}{3p}\right\rfloor.
$$

\endproof

We combine all the divisibility properties of $\prod_{i=n}^1 C_i$ together
and get the following

\begin{proposition}\label{prop1}
Let the prime factorisation of $a^*$ be $a^* =
2^{\sigma_0}p_1^{\sigma_1}p_2^{\sigma_2}\cdots p_d^{\sigma_d}$ where
$\sigma_0$ can be equal to zero while the other powers $\sigma_i$
are strictly positive. Define $\PP_1:= \{p_1, \ldots, p_d\}$,
$\PP_2:= \PP\setminus (\PP_1\cup\{2,3\})$. If $2\mid a^*$ then
\begin{equation}\label{prop1_eq1}
\gcd(p_n,q_n) \ge 2^{\big\lfloor \frac{n}{4}\big\rfloor}\prod_{p_i\in\PP_1} p_i^{\;\sum\limits_{j=1}^{\sigma_i} \left\lfloor
\frac{2n+p^j-1}{2p^j}\right\rfloor} \cdot \prod_{p\in\PP_2} p^{\big\lfloor \frac{3n+p-2}{3p}\big\rfloor}.
\end{equation}
If $2\nmid a^*$ then
\begin{equation}\label{prop1_eq2}
\gcd(p_n,q_n) \ge 2^{\big\lfloor \frac{n+3}{8}\big\rfloor}\prod_{p_i\in\PP_1} p_i^{\;\sum\limits_{j=1}^{\sigma_i} \left\lfloor
\frac{2n+p^j-1}{2p^j}\right\rfloor} \cdot \prod_{p\in\PP_2} p^{\big\lfloor \frac{3n+p-2}{3p}\big\rfloor}.
\end{equation}
\end{proposition}

We now provide shorter lower bounds for~\eqref{prop1_eq1}
and~\eqref{prop1_eq2} and then provide slightly better ones that,
after some efforts, can still be made effective for large enough
$n$. Observe that $\big\lfloor \frac{2n+ p^j-1}{2p^j}\big\rfloor \ge
\big\lfloor \frac{n}{p^j}\big\rfloor$ and $\big\lfloor \frac{3n+
p-2}{3p}\big\rfloor \ge \big\lfloor \frac{n}{p}\big\rfloor$. For
convenience, if $3\nmid a^*$ we still add $3$ to the set $\PP_1$ by
setting $p_{d+1}:=3, \sigma_{d+1}:=0$. Then
$$
2^{\sum_{i=1}^\infty \frac{n}{2^i}}\cdot \prod_{p_j\in\PP_1}
p_j^{\sum_{i=1}^{\sigma_j}\big\lfloor \frac{n}{p_j^i}\big\rfloor}
\cdot \prod_{p\in\PP_2} p^{\big\lfloor\frac{n}{p}\big\rfloor}\cdot
\prod_{p_j\in\PP_1} p_j^{\sum_{i=\sigma_j+1}^\infty \frac{n}{p_j^i}}
\cdot \prod_{p\in\PP_2}p^{\sum_{i=2}^\infty \frac{n}{p^i}} \ge n!\ge
\sqrt{2\pi n} \left(\frac{n}{e}\right)^n.
$$
The last inequality infers that
\begin{equation}\label{gcdpq}
\gcd(p_n,q_n)\ge \sqrt{2\pi n}(c_1n)^n
\end{equation}
where $c_1=c_1(a^*)$ is defined as
\begin{equation}\label{def_c1}
c_1 = \left\{ \begin{array}{ll} 2^{-3/4}
\exp\left(-1-\sum_{p_j\in\PP_1} \frac{\ln
p_j}{p_j^{\sigma_j}(p_j-1)} - \sum_{p\in \PP_2} \frac{\ln
p}{p(p-1)}\right)& \mbox{if }\; 2\mid a^*;\\[3ex]
2^{-7/8} \exp\left(-1-\sum_{p_j\in\PP_1} \frac{\ln
p_j}{p_j^{\sigma_j}(p_j-1)} - \sum_{p\in \PP_2} \frac{\ln
p}{p(p-1)}\right)& \mbox{if }\; 2\nmid a^*.
\end{array}\right.
\end{equation}
$c_1$ reaches its minimal value in the case $\PP_1 = \{3\}$ with
$\sigma_1 = 0$. Then $c_1\approx 0.0924$. However, if $a^* = 3$ then
$c_1\approx 0.1333$. In general, more squares of small prime numbers
divide $a^*$, bigger is the value of $c_1$.

We can provide a better asymptotic lower estimate on $\gcd(p_n,q_n)$
for large enough $n$. The exact condition on $n$ can be effectively
computed, however the computations will not be nice. Consider a
prime $p\in\PP_2$. The term $p^{\big\lfloor
\frac{3n+p-2}{3p}\big\rfloor}$ has an extra power of $p$ compared to
$p^{\lfloor n/p\rfloor}$ if for some integer $k$,
$$
\frac{n}{p}<k\le \frac{3n+p-2}{3p}\qquad \Longleftrightarrow\qquad
\frac{n}{k}< p\le \frac{3n-2}{3k-1}.
$$
We also have $\prod_{p\in\PP_1} p\asymp 1$ where the implied
constants only depend on $a^*$ but not on $n$. Define the set
$$
K:= \bigcup_{k=1}^n \left(\frac{n}{k}, \frac{3n-2}{3k-1}\right]
$$
Then $\gcd(p_n,q_n)\ge T \cdot \sqrt{2\pi n}(c_1n)^n$ where
$$
T \asymp \prod_{p\in \PP\cap K} p = \exp\left(\sum_{p\in K\cap\PP}
\ln p\right) = \exp \left( \sum_{k=1}^n
\left(\theta\left(\frac{3n-2}{3k-1}\right) -
\theta\left(\frac{n}{k}\right)\right)\right),
$$
where $\theta(x)$ is the first Chebyshev function. It is well known
(see~\cite{ros_sch_62} for example) that for large enough $x$,
$|\theta(x) - x|< \frac{x}{2\ln x}$. Therefore for $y>x$ one has
$\theta(y)-\theta(x)\ge y-x - \frac{y}{\ln y}$. This implies
$$
\sum_{k=1}^{\sqrt{\ln n}} \left(\theta\left(\frac{3n-2}{3k-1}\right)
- \theta\left(\frac{n}{k}\right)\right) \ge \sum_{k=1}^{\sqrt{\ln
n}} \frac{n-2k}{k(3k-1)}\ -\ O\left(\frac{n}{\sqrt{\ln n}}\right) =
n\sum_{k=1}^{\sqrt{\ln n}} \frac{1}{k(3k-1)}\ -\
O\left(\frac{n}{\sqrt{\ln n}}\right).
$$
For any $\epsilon>0$ and for large enough $n$, the last expression
can be made bigger that $(\tau - \epsilon)n$ where $\tau:=
\sum_{k=1}^\infty \frac{1}{k(3k-1)}\approx 0.74102$. Therefore $T
\gg e^{(\tau -\epsilon)n} = \gamma^{n(1-\epsilon)}$ where $\gamma =
e^\tau$. Finally, we get
\begin{equation}\label{def_c2}
\gcd(p_n,q_n)\gg ((c_2-\delta) n)^n,\quad \mbox{where } c_2 =
c_1\cdot \gamma,
\end{equation}
$\delta$ can be made arbitrarily small and the implied constant in
the inequality only depends on $a^*$ and $\delta$ but not on $n$.
For the case $a^* = 1$, when the constant $c_1$ is minimal possible,
we get $c_2 \approx 0.1939$. Respectively, for $a^*=6$, $c_2 \approx
0.2797$.

\section{Lower and upper bounds on the denominators $q_n$.}

In this section we will get upper and lower bounds of the
denominators $q_n$, compared to $q_{n-1}$. Since the recurrent
formulae between $q_n, q_{n-1}$ and $q_{n-2}$ depend on $n$ modulo
4, it makes sense to compare $q_{4k}$ and $q_{4k+4}$.

We adapt some notation from~\cite{badziahin_2022}. Denote
$$
T_{4k}:= \left(\begin{array}{cc} p_{4k}&q_{4k}\\
p_{4k-4}&q_{4k-4}
\end{array}\right)
$$
Then~\cite[(69) and (70)]{badziahin_2022} one has
\begin{equation}\label{eq2}
T_{4k+4} = \left(\begin{array}{cc}
a_{k11}&a_{k12}\\
1&0
\end{array}
\right) S_{4k}
\end{equation}
where $a_{k11}$ and $a_{k12}$ are the corresponding indices of
$C_{4k+4}C_{4k+3}C_{4k+2}C_{4k+1}$. In view of~\eqref{no4}, one
computes
\begin{equation}\label{eq3}
\begin{array}{rl}
a_{k11} = 2(8k+3)(8k+5)(8k+7)(8k+9)(t_1t_2)^2 &+ 6(8k+5)(8k+7)(36k^2+55k+16)a^*t_1t_2\\
&+ (12k+5)(12k+11)(3k+2)(6k+7)a^{*2}
\end{array};
\end{equation}$$
a_{k12} = 2(12k+1)(3k+1)(8k+7)a^*t_1((8k+5)(8k+9)t_1t_2 +
2(36k^2+63k+25)a^*).
$$
To make the notation shorter, we write $a_{k12} =
2(12k+1)(3k+1)(8k+7)a^*t_1 p(k)$ where $p(k)$ is a polynomial of $k$
with parameters $t_1t_2$ and $a^*$.

Then an easy adaptation of the proof of~\cite[Lemma
16]{badziahin_2022} gives
\begin{lemma}\label{lem5}
Let $a^*\in\NN$ and $t_1,t_2\in\ZZ$ satisfy $12a^*\le |t_1t_2|$.
Then $q_{4k+4}$ and $q_{4k}$ satisfy the relation
\begin{equation}\label{lem5_eq1}
|q_{4k+4}| > (8k+3)(8k+5)(8k+7)(8k+9) (t_1t_2+2a^*)^2 |q_{4k}|.
\end{equation}
\end{lemma}

Now we will provide an opposite inequality between the denominators
$q_{4k+4}$ and $q_{4k}$. Three consecutive denominators of this form
are related by the equation~\cite[(72)]{badziahin_2022}:
\begin{equation}\label{eq4}
q_{4k+4} = a_{k11}q_{4k} +
(dq_{4k-4}-b_{k21}q_{4k})\frac{a_{k12}}{b_{k22}},
\end{equation}
where $b_{k21}/d$ and $b_{k22}/d$ are the corresponding entries of
$C_{4k-2}^{-1}C_{4k-1}^{-1}C_{4k}^{-1}$, i.e.
\begin{equation}\label{def_d}
d = -(12k-7)(12k-5)(12k-1)(3k-1)(6k-1)(6k+1)a^{*3},
\end{equation}
$$
b_{k21} = -(12k-5)(6k-1)a^* - 2(8k-3)(8k-1)t_1t_2,
$$$$
b_{k22} = 2(8k-1)t_1((8k-3)(8k+1)t_1t_2+2(36k^2-9k-2)a^*) =:
2(8k-1)t_1p(k-1).
$$

\begin{lemma}\label{lem6}
Let $a^*$, $t_1, t_2$ be the same as in Lemma~\ref{lem5}. Then
$q_{4k+4}$ and $q_{4k}$ satisfy the following relations:
\begin{equation}\label{lem6_eq1}
|q_{4k+4}|\le 2(8k+3)(8k+5)(8k+7)(8k+9)\left(
t_1t_2+\frac{135}{32}a^*\right)^2|q_{4k}|,\quad\mbox{if }\;
t_1t_2>0,
\end{equation}
\begin{equation}\label{lem6_eq2}
|q_{4k+4}|\le 2(8k+3)(8k+5)(8k+7)(8k+9)\left(t_1t_2 + \frac{27}{32}
a^*\right)^2|q_{4k}|,\quad\mbox{if }\; t_1t_2<0.
\end{equation}
\end{lemma}

\proof First, we estimate the terms in~\eqref{eq4}. Since
$|t_1t_2|\ge 12a^*$, we get for all $k\ge 1$ that
$\frac1{12}(12k-5)(6k-1)(12a^*) < (8k-3)(8k-1)(12a^*) \le
(8k-3)(8k-1)|t_1t_2|$. Therefore
\begin{equation}\label{eq6}
|b_{k21}|\le 3(8k-3)(8k-1)|t_1t_2|.
\end{equation}
Next, by Lemma~\ref{lem5} we have
$$
|dq_{4k-4}|\le
\frac{(12k-7)(12k-5)(12k-1)(3k-1)(6k-1)(6k+1)a^{*3}}{(8k-5)(8k-3)(8k-1)(8k+1)(t_1t_2+2a^*)^2}
q_{4k}.
$$
Since $|t_1t_2 + 2a^*|\ge 10a^*$ and $12a^* \le |t_1t_2|$, one can
verify that
\begin{equation}\label{eq5}
|dq_{4k-4}| < (8k-3)(8k-1) |t_1t_2 q_{4k}|.
\end{equation}
Next, we have $(8k+5)(8k+9)|t_1t_2| > 12(36k^2+63k+25)a^*$,
therefore we always have
$$
\frac{|p(k)|}{(8k+5)(8k+9)|t_1t_2|}\in \left\{\begin{array}{ll}
\left[ 1,\frac76 \right]&\mbox{if }\; t_1t_2\ge0;\\[1ex]
\left[\frac56, 1\right]& \mbox{if }\; t_1t_2<0.
\end{array}\right.
$$
The last inequality in turn implies that for $k\ge 1$ the ratio
$a_{k12}/b_{k22}$ is always positive and satisfies
\begin{equation}\label{eq10}
\frac{a_{k12}}{b_{k12}} \le \frac{6(12k+1)(3k+1)(8k+7)a^*}{8k-1}.
\end{equation}

Assume that $t_1t_2\ge 0$. In that case, the last inequality
together with~\eqref{eq6} and~\eqref{eq5} imply that
$$
\left| (dq_{4k-4}-b_{k21}q_{4k})\frac{a_{k12}}{b_{k22}}\right| \le
24 (8k-3)(8k+7)(12k+1)(3k+1) a^*t_1t_2 q_{4k}.
$$
One can check that for all $k\ge 1$,
\begin{equation}\label{eq7}
\frac{6(8k+5)(8k+7)(36k^2+55k+16) + 24
(8k-3)(8k+7)(12k+1)(3k+1)}{2(8k+3)(8k+5)(8k+7)(8k+9)}
<\frac{135}{16}
\end{equation}
and
\begin{equation}\label{eq8}
\frac{81}{256}<\frac{(12k+5)(12k+11)(3k+2)(6k+7)}{2(8k+3)(8k+5)(8k+7)(8k+9)} \le
\frac{23}{66} < 1.
\end{equation}
These bounds together with the formula~\eqref{eq3} and
equation~\eqref{eq4} imply the inequality~\eqref{lem6_eq1} for $k\ge
1$. Finally, this bound can be easily verified for $k=0$ from the
equation $q_4 = a_{011}q_0 + a_{012} q_{-1}$ and $q_{-1} = 0$.

Consider the case $t_1<0$. One can check that for all $k\ge 1$,
\begin{equation}\label{eq9}
\frac{321}{187}\ge
\frac{6(8k+5)(8k+7)(36k^2+55k+16)}{2(8k+3)(8k+5)(8k+7)(8k+9)} \ge
\frac{27}{16}.
\end{equation}
This together with the condition $|t_1t_2|>12a^{*2}$ imply that
$a_{k11}>0$ and $q_{4k}$ and $q_{4k+4}$ share the same sign for all
$k\in\NN$. Next, since $(12k-5)(6k-1)a^* < (8k-3)(8k-1)|t_1t_2|$, we
have that $b_{k21}>0$ and then in view of~\eqref{eq5} and
$\frac{a_{k12}}{b_{k22}}>0$, the term
$(dq_{4k-4}-b_{k21}q_{4k})\frac{a_{k12}}{b_{k22}}$ has the opposite
sign compared to $a_{k11}q_{4k}$. That all implies that
$|q_{4k+4}|\le |a_{k11} q_{4k}|$. Finally, the
inequalities~\eqref{eq8} together with~\eqref{eq9} establish the
bound~\eqref{lem6_eq2}.\endproof

Lemma~\ref{lem6} immediately implies that for $t_1t_2>0$,
$$
|q_{4k}|\le 2^k \left(t_1t_2 +
\frac{135}{32}a^*\right)^{2k}(8k+1)!! \le 16k \left(8\cdot
2^{1/4}e^{-1}\sqrt{t_1t_2 + \frac{135}{32}a^*}\right)^{4k} k^{4k}.
$$
The case of $t_1t_2<0$ can be dealt with in a similar way. Finally,
we get the estimate
\begin{equation}\label{ubound_q}
|q_{4k}| \le 16 k c_3^{4k} k^{4k},
\end{equation}
where
$$
c_3 = c_3(t_1,t_2,a^*) =
\left\{
\begin{array}{ll} 8\cdot 2^{1/4}e^{-1}\sqrt{t_1t_2 +
\frac{135}{32}a^*}&\mbox{if }\;
t_1t_2>0\\[2ex]
8\cdot 2^{1/4}e^{-1}\sqrt{|t_1t_2| - \frac{27}{32}a^*}&\mbox{if }\;
t_1t_2<0.
\end{array}\right.
$$

\begin{lemma}\label{lem7}
Under the same conditions on $a^*, t_1,t_2$ as in the previous lemma, one has
\begin{equation}\label{lem7_eq1}
|q_{4k+4}|\ge 2(8k+3)(8k+5)(8k+7)(8k+9)\left(
t_1t_2+\frac{9}{16}a^*\right)^2|q_{4k}|,\quad\mbox{if }\; t_1t_2>0,
\end{equation}
\begin{equation}\label{lem7_eq2}
|q_{4k+4}|\ge 2(8k+3)(8k+5)(8k+7)(8k+9)\left(t_1t_2 + 3
a^*\right)^2|q_{4k}|,\quad\mbox{if }\; t_1t_2<0.
\end{equation}
\end{lemma}

\proof If $t_1t_2>0$ we have $q_{4k+4} > a_{k11}q_{4k}$.  Then the
lower bound in~\eqref{eq8} together with the lower bound
in~\eqref{eq9} imply the bound~\eqref{lem7_eq1}.

Now assume that $t_1t_2<0$. Then, as we have shown in the proof of
Lemma~\ref{lem6}, $b_{k21}>0$ and
$dq_{4k-4}$ and $b_{k21}q_{4k}$ have the opposite signs. This
together with $\frac{a_{k12}}{b_{k22}}>0$, the
inequality~\eqref{eq10} and $0<b_{k21}\le 2(8k-3)(8k-1)|t_1t_2|$ in
turn imply that
$$
|q_{4k+4}|\ge |a_{k11}q_{4k}| - \frac{b_{k21}a_{k12}}{b_{k22}}|q_{4k}| \ge (a_{k11} + 12(8k-3)(12k+1)(3k+1)(8k+7)a^*t_1t_2)|q_{4k}|
$$
We need to show that the expression
$$a_{k11}+
12(8k-3)(12k+1)(3k+1)(8k+7)a^*t_1t_2 - 2(8k+3)(8k+5)(8k+7)(8k+9)\left(t_1t_2
+ 3 a^*\right)^2
$$
is always positive. Notice that after substituting~\eqref{eq3} into it and
expanding the brackets, the term with $(t_1t_2)^2$ disappears. The term for
$a^*t_1t_2$ then equals to
$$
-6(8k+7)(160k^3+1532k^2+1063k+196)a^*t_1t_2
$$
and the term for $a^{*2}$ is
$$-(71136 k^4 + 212976 k^3 + 227970 k^2 + 102633 k + 16240)
$$
(we made these computations with Wolfram Mathematika). Finally, one can check
that in the case $|t_1t_2|> 12a^*$, the absolute value of the first term is
always bigger than that of the second term and therefore the whole expression
is positive.

{\bf Remark.} By performing neater computations, one can make the coefficient
3 in $(t_1t_2+3a^*)^2$ slightly smaller. However we decide not to further
complicate already tedious calculations.
\endproof

Analogously to~\eqref{ubound_q}, one can find shorter lower bounds
for $|q_{4k}|$. With help of the known inequality $(8k+1)!!\ge
8k(8k/e)^{2k}$, Lemma~\ref{lem7} infers
\begin{equation}\label{lbound_q}
|q_{4k}|\ge 8 k c_4^{4k} k^{4k},
\end{equation}
where
$$
c_4 = c_4(t_1,t_2,a^*) =
\left\{
\begin{array}{ll} 8\cdot 2^{1/4}e^{-1}\sqrt{t_1t_2 +
\frac{9}{16}a^*}&\mbox{if }\;
t_1t_2>0\\[2ex]
8\cdot 2^{1/4}e^{-1}\sqrt{|t_1t_2| - 3a^*}&\mbox{if }\; t_1t_2<0.
\end{array}\right.
$$

\section{Distance between $x$ and the convergents}

From~\cite[Lemma 17]{badziahin_2022} we know that, under the condition
$12a^*\le |t_1t_2|$, one has
$$
\left|x - \frac{p_{4k}}{q_{4k}}\right| < 2 \left|\frac{p_{4k}}{q_{4k}} -
\frac{p_{4k+4}}{q_{4k+4}}\right|.
$$
In order to estimate the right hand side, we use the matrix
equation~\cite[(72)]{badziahin_2022}:
$$
T_{k+1} = \left(\begin{array}{cc}
a_{k11} - a_{k12} \frac{b_{k21}}{b_{k22}}&\frac{da_{k12}}{b_{k22}}\\
1&0
\end{array}\right) T_k.
$$
Notice that the values of $d$ in fact depends on $k$ (see the
formula~\eqref{def_d}). to emphasize this dependence, in this
section we write $d(k)$ for it. Then the above equation gives the
following formula:
$$
\left|\frac{p_{4k}}{q_{4k}} -
\frac{p_{4k+4}}{q_{4k+4}}\right| = \frac{\left|\prod_{i=1}^k \frac{d(i)a_{i12}}{b_{i22}}\right|\cdot |p_0q_4-q_0p_4|}{q_{4k}q_{4k+4}}.
$$
We first compute its product term: {\footnotesize$$
\left|\prod_{i=1}^k \frac{d(i)a_{i12}}{b_{i22}}\right| =
\prod_{i=1}^k
\frac{2(12i-7)(12i-5)(12i-1)(3i-1)(6i-1)(6i+1)(12i+1)(3i+1)(8i+7)a^{*4}|t_1p(i)|}{2(8i-1)|t_1p(i-1)|}
$$}
$$
= \frac{(8k+7)|p(k)|}{7|p(0)|}\cdot
\frac{(3k+1)(6k+1)(12k+1)a^{*4k}(12k)!}{2^{6k}3^{4k}(4k)!} $$$$\le
\frac{\sqrt{3}(3k+1)(6k+1)(12k+1)(8k+7)|p(k)|}{7|p(0)|}\cdot
\left(\frac{12^2k^2a^*}{2\sqrt{2}e^2}\right)^{4k}.
$$
Next, from~\eqref{eq2} for $k=0$ we get that $|p_0q_4-p_4q_0| =
14a^*t_1|p(0)|$. Finally, we unite all these bounds together with
the lower bounds~\eqref{lem7_eq1},~\eqref{lem7_eq2}
and~\eqref{lbound_q} for $|q_{4k}|$ to get $$ \left|x -
\frac{p_{4k}}{q_{4k}}\right| \le
\frac{2\sqrt{3}(3k+1)(6k+1)(12k+1)(8k+7)|t_1a^*p(k)|}{2(8k+3)(8k+5)(8k+7)(8k+9)(|t_1t_2|-3a^*)^2\cdot
64k^2}\cdot\left(\frac{12^2a^*}{2\sqrt{2}e^2c_4^2}\right)^{4k}
$$

To simplify the right hand side, notice that
$\frac{(3k+1)(6k+1)(12k+1)}{(8k+3)(8k+5)(8k+9)} < \frac{27}{64}$.
Next, since $|t_1t_2|\ge 12a^*$, one has $|p(k)| =
|(8k+5)(8k+9)t_1t_2 + 2(36k^2+63k+25)a^*|\le 2(8k+5)(8k+9)|t_1t_2|$
which for all $k\ge 1$ is smaller than $442k^2|t_1t_2|$. Finally,
$(|t_1t_2|-3a^*)^2 \ge \frac9{16} (t_1t_2)^2$. Collecting all of
these inequalities together gives,
\begin{equation}\label{ubound_dist}
\left|x - \frac{p_{4k}}{q_{4k}}\right| \le \frac{\sqrt{3}\cdot
27\cdot 442|t_1^2t_2a^*|}{64\cdot (9/16)\cdot 64 (t_1t_2)^2}\cdot
\left(\frac{72a^*}{\sqrt{2}e^2c_4^2}\right)^{4k} \le |t_1| c_5^{4k},
\end{equation}
where
$$
c_5 = \left\{
\begin{array}{ll} \frac{9a^*}{(16t_1t_2 +
9a^*)}&\mbox{if }\;
t_1t_2>0\\[2ex]
\frac{9a^*}{16(|t_1t_2| - 3a^*)}&\mbox{if }\; t_1t_2<0.
\end{array}\right.
$$

\section{Estimating the irrationality exponent}

In this section we establish Theorems~ref{th1} and~\ref{th2}.
Consider $p_k^*:=p_{4k}/\gcd(p_{4k}, q_{4k})$ and
$q_k^*:=q_{4k}/\gcd(p_{4k}, q_{4k})$. Definitely, they are both
integers and~\eqref{ubound_q} together with~\eqref{gcdpq} imply
\begin{equation}\label{bnd_q}
|q^*_k|\le 4\sqrt{\frac{2k}{\pi}} \left(\frac{c_3}{4c_1}\right)^{4k}
=: 4\sqrt{\frac{2k}{\pi}}\cdot  c_6^{4k}.
\end{equation}
For arbitrary $\delta>0$ and large enough $k$, one can use the
inequality~\eqref{def_c2} to get
\begin{equation}\label{bnd_qstar}
|q^*_k|\ll 16k\left(\frac{c_3}{4(c_2-\delta)}\right)^{4k} \ll
\left(\frac{c_3}{4c_2} + \delta_1\right)^{4k} =:
(c^*_6+\delta_1)^{4k}
\end{equation}
where $\delta_1>0$ can be made arbitrarily close to zero for large enough
$k$. Denote the upper bound for $b^*_k$ by $Q(k,t,a)$.

Next, we combine the last two inequalities with~\eqref{ubound_dist} and get
\begin{equation}\label{bnd_dx}
||q_k^*x||\le |t_1| c_5^{4k}\cdot
4\sqrt{\frac{2k}{\pi}}\left(\frac{c_3}{4c_1}\right)^{4k} \le
4|t_1|\sqrt{k}\left(\frac{c_3c_5}{4c_1}\right)^{4k} =:
4|t_1|\sqrt{k} c_7^{-4k}
\end{equation}
or
\begin{equation}\label{bnd_dxstar}
||q_k^*x||\ll \left(\frac{4c_2}{c_3c_5}-\delta_2\right)^{-4k} =:
(c_7^* - \delta_2)^{-4k}
\end{equation}
where $\delta_2$ can be made arbitrarily small and $k$ is large
enough, depending on $\delta_2$. Denote the upper bound of
$||q^*_kx||$ by $R(k,t,a)$.

Consider an arbitrary $q \ge \frac{1}{2R(1,t,a)}= q_0$. We now
impose the condition $c_7 > e^{1/4}$. In this case, by examining the
derivative of $\sqrt{k}c_7^{-4k}$, one can check that it strictly
decreases for $k\ge 1$. Therefore, there exists a unique $k\ge 2$
such that $R(k,t,a)< \frac{1}{2q}\le R(k-1,t,a)$. Let $p\in \ZZ$ be
such that $||qx|| = |qx-p|$. Since two vectors $(p_k^*, q_k^*)$ and
$(p_{k+1}^*, q_{k+1}^*)$ are linearly independent, at least one of
them must be linearly independent with $(p,q)$. Suppose that is
$(p_k^*, q_k^*)$. Then we estimate the absolute value of the
following determinant:
$$
1\le \left| \begin{array}{cc}
q&q_k^*\\
p&p_k^*
\end{array}\right| \le \left| \begin{array}{cc}
q&q_k^*\\
p-qx& p_k^* - q_k^*x
\end{array}\right| \le q R(k,t,a) + ||qx||Q(k,t,a).
$$
Since $qR(k,t,a) < \frac12$, we must have $||qx|| \ge (2Q(k,t,a))^{-1}$.
Analogously, if $(p,q)$ is linearly independent with $(p_{k+1}^*,
q_{k+1}^*)$, we have $||qx||\ge (2Q(k+1,t,a))^{-1}$. The latter lower bound
is weaker. Now, we need to rewrite the right hand side of the inequality in
terms of $q$ rather than $k$.

Since $\frac{1}{2q}\le R(k-1,t,a)$, we have that
$$
\frac{c_7^{4(k-1)}}{8|t_1| \sqrt{k-1}} \le q\quad
\Longrightarrow\quad k-1\le \frac{\log
(8|t_1|q)+\log\log(8|t_1|q)}{4\log c_7}.
$$
The last implication can be justified by standard techniques on
working with logarithms, see~\cite[(41)]{badziahin_2022}.

Finally, substitute the last lower bound for $k$ in $||qx||\ge
(2Q(k+1,t,a))^{-1}$ and get
$$
||qx|| \ge \frac{\sqrt{\pi}}{8\sqrt{2(k+1)} c_6^{4(k+1)}} \ge
\frac{2\sqrt{\pi} (\log c_7)^{1/2}}{8\sqrt{6}c_6^8 (\log (8|t_1| q)
+ \log\log(8|t_1|q))^{1/2}\cdot (8|t_1|q)^{\frac{\log c_6}{\log
c_7}} (\log (8|t_1| q))^{\frac{\log c_6}{\log c_7}}}.
$$$$
\ge \frac{(\log c_7)^{1/2}}{8c_6^8(8|t_1|)^{\frac{\log c_6}{\log
c_7}}}\cdot q^{-\frac{\log c_6}{\log c_7}} (\log (8|t_1|
q))^{-\frac{\log c_6}{\log c_7} - \frac12} =
\frac{\tau}{g_2}q^{-\lambda}(\log(8|t_1|q))^{-\lambda-\frac12}.
$$
To finish the proof of Theorem~\ref{th1}, we recall, that for
convenience, we in fact worked with the number $x/g_2$ rather than
$x$, i.e. the inequality above is for $||qx/g_2||$. Hence one needs
to multiply both sides by $g_2$.

Regarding theorem~\ref{th2}, we use inequalities~\eqref{bnd_qstar}
and~\eqref{bnd_dxstar}. in this case the computations are much
easier and we get for any $\delta_3>0$ and large enough integer $q$
that
$$
||qx|| \ge q^{-\frac{\log c_6^*}{\log c_7^*} - \delta_3}
$$
or in other words $\lambda_{eff}(x)\le \frac{\log c_6^*}{\log
c_7^*}$. That completes the proof of Theorem~\ref{th2}.

\end{document}